\documentclass[a4paper]{article}
\usepackage[english]{babel}
\usepackage[utf8]{inputenc}
\usepackage{amsmath,amsthm}
\usepackage{amsfonts,amssymb}

\usepackage{tikz,tikz-cd} 

\usepackage{xifthen} 

\newcommand{\toposdefaut}{0}
\newcommand{\topos}[1][\toposdefaut]{ 
\ifthenelse{\equal{#1}{0}}{ \mathcal{T} }
{
\ifthenelse{\equal{#1}{1}}{ \mathcal{E} }{ #1 }
}
}

\newcommand{\sh}{\textsf{Sh}}

\newcommand{\dom}{\text{Dom}}

\newcommand{\spec}{\text{Spec }}
\newcommand{\spect}{\text{Spec}^{\infty}}

\newcommand{\Q}{\mathbb{Q}}
\newcommand{\R}{\mathbb{R}}
\newcommand{\C}{\mathbb{C}}
\newcommand{\N}{\mathbb{N}}

\newcommand{\Ecal}{\mathcal{E}} 
 
\newcommand{\Tcal}{\mathcal{T}}

\newcommand{\Ocal}{\mathcal{O}}

\newcommand{\Lcal}{\mathcal{L}}

\newcommand{\Ccal}{\mathcal{C}}

\usepackage{titlesec}
\titleformat{\subsection}[runin]{\normalfont}{\thesubsection}{0pt}{}[.]

\newcommand{\block}[1]
{

\par \subsection{} #1

\bigskip}

\newcommand{\blockn}[1]{\par #1 \bigskip}

\newcommand{\Th}[1]
	{
	\bigskip	
	\textbf{Theorem : }{\itshape #1}
		
	\bigskip
	}

\newcommand{\Prop}[1]
	{

	\bigskip
	
	\textbf{Proposition : }{\itshape #1}
		
	\bigskip
	
	}

\newcommand{\Cor}[1]
	{

	\bigskip
	
	\textbf{Corollary : }{\itshape #1}	
		
	\bigskip

	}

\newcommand{\Lem}[1]
	{

	\bigskip
	
	\textbf{Lemma : }{\itshape #1}
		
	\bigskip
	
	}

\newcommand{\Def}[1]
	{
	
	\bigskip
	
	\textbf{Definition : }{\itshape #1}
	
	\bigskip
	
	}

\newcommand{\Dem}[1]{
	
	\smallskip
	
	\textbf{Proof : } \par
	 {#1} $\square$
	 
	 \bigskip
}

\begin{document}

\pagestyle{plain}
\title{Constructive Gelfand duality for non-unital commutative $C^{*}$-algebras}
\author{Simon Henry}

\renewcommand{\thefootnote}{\fnsymbol{footnote}} 
\footnotetext{\emph{email:} simon.henry@imj-prg.fr}
\renewcommand{\thefootnote}{\arabic{footnote}} 


\maketitle

\begin{abstract}
We prove constructive versions of various usual results related to the Gelfand duality. Namely, that the constructive Gelfand duality extend to a duality between commutative nonunital $C^{*}$-algebras and locally compact completely regular locales, that ideals of a commutative $C^{*}$-algebras are in order preserving bijection with the open sublocales of its spectrum, and a purely constructive result saying that a commutative $C^{*}$-algebra has a continuous norm if and only its spectrum is open. We also extend all these results to the case of localic $C^{*}$-algebras. In order to do so we develop the notion of one point compactification of a locally compact regular locale and of unitarization of a $C^{*}$-algebra in a constructive framework.
\end{abstract}

\tableofcontents

\section{Notations and preliminares}

\blockn{This paper has been written to provide two technical tools which were needed in the proof of the main theorems of \cite{henry2015toward} : the non-unital Gelfand duality, including the characterization of the spectrum given by proposition \ref{spec_class_char}, and (one direction) of the ``possitivity" theorem  \ref{continuity=openess}. We took the opportunity to prove (constructively) some other results in this spirit that might be useful for future works, like for example theorem \ref{openeqideal} and the results of section \ref{secLocalic}. }

\blockn{In all this paper we are working in the internal logic of an elementary topos with a natural number object $\N$. The subobject classifier is denoted by $\Omega$, and $\top$ and $\bot$ denotes its top and bottom element, i.e. the proposition true and false.}

\blockn{A frame is a complete Heyting algebra, a frame homomorphism is an order preserving map commuting to arbitrary supremums and finite infimums. The category of locales is defined as the opposite of the category of frames, if $X$ is a locale the corresponding frame is denoted by $\Ocal(X)$. If $f$ is a morphism of locales, the corresponding frame homomorphism is denoted $f^{*}$.

Elements of $\Ocal(X)$ are called open sublocales of $X$. The top element of $\Ocal(X)$ is denoted by $X$, the bottom element by $\emptyset$. When talking about open sublocales, ``$V$ is bigger than $U$" or ``$U$ is smaller than $V$" always means $U \leqslant V$. For more information on the theory of locale, the reader can consult \cite{picado2012frames} (which is unfortunately non constructive) or \cite[C1]{sketches}. 

Supremums and finite infimums in $\Ocal(X)$ are called unions and intersections and are denoted by the symbols $\cup$ and $\cap$.
}

\blockn{If $U$ is an open sublocale of $X$, then $\neg U$ denote the open sublocale $U \Rightarrow \emptyset$ and $U^{c}$ denote the closed complement of $U$, ie the locale such that $\Ocal(U^{c}) = \{ V \in \Ocal(X) | U \leqslant V \}$. In particular, $\neg U$ is the interior of $U^{c}$.}

\blockn{An (increasing) net of open sublocales of $X$ is an \emph{inhabited} family $(U_i)_{i \in I}$ of open sublocales of $X$ such that for each $i,j \in I$ there exists $k \in I$ such that $U_k$ is bigger than $U_i$ and $U_j$.}

\blockn{If $U$ and $V$ are two open sublocales of a locale $X$, we say that: 

\begin{itemize}
\item $U \ll V $ ($U$ is way below $V$) if for each increasing net $(U_i)_{i \in I}$ whose supremum is bigger than $V$ there exists $i \in I$ such that $U \leqslant U_i$.

\item $ U \triangleleft V$ ($U$ is rather below $V$) if there exists $W \in \Ocal(X)$ such that $V \cup W =X$ and $U \cap W = \emptyset$. Or equivalently if $\neg U \cup V = X$.

\item $U \triangleleft_{CR} V$ ($U$ is completely below $V$) if there exists a ``scale" $(U_q)_{q \in [0,1] \cap \Q}$ such that $U_0=U$, $U_1=V$ and for each $q'<q$ one has $U_{q'}  \triangleleft U_q$. This is also equivalent to the existence of a function $f$ from $X$ to the locale\footnote{We mean the formal locale of real number, which might be non spatial and hence different from the topological space of real number in the absence of the law of excluded middle.} $[0,1]$ of real numbers between $0$ and $1$ such that $f^{*}( ]0,1]) \subset V$ (i.e. $f$ restrictd to $V^{c}$ is $0$) and $f$ restricted to $U$ is constant equal to $1$. (see \cite[V.5.7 and XIV.6.2]{picado2012frames}).

\end{itemize}

One also says that $X$ is locally compact (resp. regular, resp. completely regular) if any open sublocale $V$ of $X$ can be written as a supremum of open sublocale $U$ such that $U \ll V$ (resp. $U \triangleleft V$, resp. $U \triangleleft_{CR} V$). A locale $X$ is said to be compact if $X \ll X$.

\blockn{One has the following properties:

\begin{enumerate}

\item Each of the three relations $\ll$, $\triangleleft$ and $\triangleleft_{CR}$ satisfies the properties: if $a\leqslant b$, $b \ll c$ and $c \leqslant d$ then  $a \ll d$; and if $a \ll b$ and $c \ll d$ then $a \cup c \ll b \cup d$.

\item In a regular (resp. completely regular) locale $U \ll V$ implies $U \triangleleft V$ (resp. $U \triangleleft_{CR} V$ ).

\item In a locally compact locale $X$, if $a \ll b$ then there exists $c$ such that $a \ll c \ll b$.

\item In a compact locale $a \triangleleft b$ imply $a \ll b$, more generally in any locale $X$, if $a \triangleleft b$ and $a \ll X$ then $a \ll b$.

\item In particular, A compact regular locale is locally compact.

\item And also, in a compact regular locale, $\triangleleft$ and $\ll$ are equivalent and in a compact completely regular locale these three relations are equivalent.

\end{enumerate}

}

}

\blockn{If $f:\Tcal \rightarrow \Ecal$ is a geometric morphism between two toposes, and $X$ a locale in the internal logic of $\Ecal$, we denote by $f^{\sharp}(X)$ the pullback of $X$ along $f$, in particular $f^{*}(\Ocal(X))$ is different from $\Ocal(f^{\sharp}(X))$ but is still a basis\footnote{This means that $f^{*}(\Ocal(X))$ generate $\Ocal(f^{\sharp}(X))$ under arbitrary join} of the topologies of $f^{\sharp}(X)$. }

\blockn{We will frequently use expression of the form $\bigcup_{u} a$ where $u$ is a proposition and $a$ is an element of a frame which might seems strange for a reader unfamiliar with this. This expression makes sense because, as a proposition, $u$ is a subset of the singleton and $a$ can be seen as a family of elements indexed by the singleton (and hence also by its subset $u$). This is of course the same as $a \cap p^{*}(u)$ where $p^{*}$ denotes the canonical frame homomorphism from the initial frame $\Omega$. But the expression with a union allows to emphasize the fact that this is indeed a union, and also it might happen that the expression defining ``$a$" only makes sense when $u$ holds, in which case only the first expression makes sense.}

\blockn{We conclude these preliminaries by the definition of real numbers. We will need two spaces of real numbers.

The first one is the set of non-negative upper semi-continuous real numbers, where the norm function of $C^{*}$-algebras will take value. A upper-semi-continuous real number is a subset $x$ of the set $\Q$ of rational such that:

\begin{itemize}
\item $\exists q \in x$
\item If $q \in x$ and $q<q'$ then $q' \in x$
\item For all $q \in x $ there exists a $q'<q$ such that $q' \in x$.
\end{itemize}

It is said to be non-negative if it is included in the set $\Q_+^{*}$ of positive rational numbers. Of course $q \in x$ has to be interpreted as $x<q$ (and will be denoted this way). Upper semi-continuous real number have good order property (every bounded set has a supremum) but poor algebraic property: even if one removes the positivity assumption, there is no opposite of an arbitrary element (the opposite of a upper semi-continuous real numbers should be a lower semi-continuous numbers) and we can only multiply positive elements.

The second is the set $\R$ of continuous real numbers, which will play the role of the ``scalars field" for $C^{*}$-algebras. A continuous real number is a pair $x=(L,U)$ of subsets of the set of $\Q$ of rational numbers such that:

\begin{itemize}

\item $U$ is upper semi-continuous real number, and $L$ is a lower semi-continuous number (i.e. satisfy the same three axioms but for the reverse order relation).

\item $L \cap U = \emptyset$.

\item for all $q<q'$ either $q \in L$ or $q' \in U$.

\end{itemize}

Of course $q \in L$ mean $q<x$ and $q \in U$ means $x<q$. The continuous real number have good algebraic properties (they form a locale ring) and topological properties (they are complete, in fact they are exactly the completion of $\Q$ by Cauchy filters) but no longer have supremums in general. The complex numbers are defined as $\R \times \R$ endowed with their usual product.

Finally the map $(L,U) \mapsto U$ induce an injection of the continuous real numbers into the semi-continuous real numbers, in particular it makes sense to wonder whether a given semi-continuous real number is continuous or not.
}

\section{One point compactification of  locales}

\blockn{In this section we will define a constructive and pointfree version of the process of one point compactification of a locally compact separated topological space.}

\block{Let $X$ be a locally compact regular locale, and $U \in \Ocal(X)$. We will denote by $\omega(U)$ the proposition:

\[ \omega(U) := ``\exists W \in \Ocal(X) \text{ such that }  W \ll X \text{ and } U \cup W =X " \]

i.e. $\omega(U)$ is the proposition ``$U$ has a compact complement". The underlying idea is that in the one point compactification, the neighbourhoods of $\infty$ are exactly the open subspaces whose complement is compact, i.e. the $U$ such that $\omega(U)$.
}

\blockn{The main result of this section is:}
\block{\label{OnepointCptmain}\Th{Let $X$ be a locally compact regular locale. Then there exists a unique (up to unique isomorphism) compact regular locale $X^{\infty}$ with a (closed) point $\{ \infty \} \subset X^{\infty}$ and an isomorphism between $X$ and the open complement of $\{ \infty \}$.

Moreover: 
\[\Ocal(X^{\infty}) \simeq \{(U,p) \in \Ocal(X) \times \Omega | p \Rightarrow \omega(U) \} \]

And the two projections from $\Ocal(X^{\infty})$ to $\Ocal(X)$ and $\Omega$ are the frame homomorphisms corresponding to the injections of $X$ and $\{\infty\}$ into $X^{\infty}$.
}

The proof will be completed in \ref{OnepointCptproof}. One can also note that this is a special case of Artin gluing\footnote{We mean by that the localic form of construction like \cite[A2.1.12, A4.1.12 and A4.5.6]{sketches}} of a closed point to $X$. This will be extremly apparent in \ref{defofXinfty} and in \ref{OnepointCptproof}.

For the rest of this section, we fix a locally compact regular locale $X$.}

\block{\Prop{The function $\omega : \Ocal(X) \rightarrow \Omega = \Ocal(\{ \infty \})$ is cartesian, i.e. order preserving and satisfies $\omega(X)=\top$ and $\omega(a \cap b) = \omega(a) \wedge \omega(b)$.}

\Dem{if $A \leqslant B$ and $\omega(A)$, then $\omega(B)$ also holds with the same $W$. One has $\omega(X)$ with $W=\emptyset$. As $\omega$ is order preserving one has $\omega(A \cap B) \leqslant \omega(A) \wedge \omega(B)$. For the converse inequality, if one has $\omega(A)$ and $\omega(B)$ then there is $W$ and $W'$ such that $W,W' \ll X$ and $W \cup A = X$ , $W \cup B =X$. Taking $W''=W \cup W'$ one has $W'' \ll X$ and $W'' \cup (A \cap B) = X$ which proves $\omega(A \cap B)$. }

}

\block{\label{defofXinfty}\Cor{There is a locale $X^{\infty}$ such that 

\[\Ocal(X^{\infty}) = \{(U,p) \in \Ocal(X) \times \Omega | p \Rightarrow \omega(U) \}\]

as in theorem \ref{OnepointCptmain}. Moreover the two projections from $\Ocal(X^{\infty})$ to $\Ocal(X)$ and $\Omega = \Ocal(\{\infty \})$ corresponds to an open inclusion of $X$ into $X^{\infty}$ and the complementary closed inclusion.  }

\Dem{From the fact that $\omega$ is cartesian one deduces\footnote{This is exactly the general construction of an Artin gluing.} that $\{(U,p) \in \Ocal(X) \times \Omega | p \Rightarrow \omega(U) \}$ is stable under arbitrary joins and finite meets in $\Ocal(X)\times \Omega$ hence it is a frame and the two projections are frame homomorphisms. Consider the element $X_0 :=(X,\bot) \in \Ocal(X^{\infty})$ then the elements of $\Ocal(X^{\infty})$ smaller than $X_0$ are the $(U,\bot)$ for $U \in \Ocal(X)$ hence the open sublocale $X_0$ is isomorphic to $X$. Conversely, the element of $\Ocal(X^{\infty})$ bigger than $X_0$ are exactly the $(X,p)$ hence the closed complement of $X_0$ is just a point, denoted $\infty$ and corresponding to $\infty^{*}(X,p)=p$.
}

For now on, the open sublocale $X_0$ of $X^{\infty}$ will be identified with $X$ (and in particular denoted $X$).

}

\block{\Lem{$X^{\infty}$ is compact.}

\Dem{Let $(X_i)_{i\in I} = (U_i,p_i)_{i \in I}$ be a covering net of open sublocales of $X^{\infty}$. As supremum in $\Ocal(X^{\infty})$ are computed componentwise one has in particular $\bigvee p_i = \top$ i.e. there exists $i_0 \in I$ such that $p_{i_0}$ holds. As $ p_i \Rightarrow \omega(U_i)$ one also has $\omega(U_{i_0})$ i.e. there exists a $W$ such that $W \ll X$ and $W \bigcup U_{i_0} =X$. As the $U_i$ form a covering net of $X$ one also has a $i_1$ such that $W \leqslant U_{i_1}$ and a $i_2$ bigger than $i_1$ and $i_0$. Hence $U_{i_2} \geqslant W \cup U_{i_0}=X$ and $p_{i_2}\geqslant p_{i_0} = \top$, hence $X_{i_2} = X^{\infty}$ which concludes the proof.
}

}

\block{\label{OnepointCptRegular}\Lem{$X^{\infty}$ is regular.}

\Dem{Let $A=(U,p)$ be any open sublocale of $X^{\infty}$. We can first see that:

\[ A = (U, \bot) \cup \bigcup_{p \atop W \ll X, U \cup W = X} (\neg W, \top) \] 

The term in the union makes sense because if $W \ll X$ then there exists a $W'$ such that $W \ll W' \ll X$ hence $W \triangleleft W'$ and $\neg W \cup W' = X$ hence $\omega(\neg W)$. $A$ is bigger than this union because when $W \cup U = X$ one has $\neg W \leqslant U$ (and when $p$ holds then $\top \leqslant p$). Conversely, using the fact that unions in $\Ocal(X^{\infty})$ are computed componentwise one easily checks that the right hand side union is indeed smaller than $(U,p) = A$.

Now, for any $V \ll U$ in $X$ one has $\omega(\neg V)= \top$ (using a $W$ such that $V \ll W \ll X$). Hence $(\neg V, \top) \in \Ocal(X^{\infty})$ is an open such that $(V, \bot) \cap (\neg V, \top)=\emptyset$ and $(\neg V, \top) \cup (U,\bot) = (X, \top)$ hence $(V, \bot) \triangleleft (U,\bot)$ and
\[ (U,\bot) = \bigcup_{V \ll U} (V,\bot) \]
because $X$ is locally compact.

Moreover, if we assume $p$, then for any $W$ such that $W \ll X$ and $(U \cup W) = \top$ one has $(\neg W,\top) \triangleleft (U, \top)=(U,p)=A$ as attested by $(W,\bot)$. Hence, if one writes: 

\[ A =  \left( \bigcup_{V \ll U} (V,\bot) \right) \cup \left( \bigcup_{p \atop W \ll X, U \cup W = \top} (\neg W, \top) \right), \]

then all the terms of the union are rather below $A$ ($\triangleleft A$) which proves that $X^{\infty}$ is regular.
}

}

\block{\label{OnepointCptproof} At this point, the existence part of theorem \ref{OnepointCptmain} and the additional properties of $X^{\infty}$ stated in \ref{OnepointCptmain} are proved, all that remains to do is to prove the uniqueness, and that is what we will do now:

\Dem{Let $Y$ be a compact regular locale with a (closed) point denoted $\infty$ such that $Y-\{\infty\}$ is identified with $X$. Let $i$ be the inclusion of $X$ into $Y$, open sublocales of $X$ will be identified with the corresponding open sublocales of $Y$ included in $X$. We will first show that for any $U \in \Ocal(X)$ one has $\infty \in i_*(U)$ if and only if $\omega(U)$ (i.e. $\omega$ is the unique ``gluing function" giving rise to a compact regular Artin gluing).

\bigskip

Indeed, assume $\omega(U)$, i.e. that there exists a $W \in \Ocal(X)$ such that $W \ll X$ and $W \cup U =X$.

Now as $W \ll X$ in $\Ocal(X)$ one also have $W \ll X$ in $Y$, hence (as $Y$ is regular) there exists $W' \in \Ocal(Y)$ such that $X \cup W' = Y$, (i.e. $\infty \in W'$) and $W \cap W' = \emptyset$. In particular $i^{*}(W') \cap W= i^{*}(W') \cap i^{*}(W) = \emptyset$ hence as $W \cup U =X$ one has $i^{*}(W') \subset U$ hence $W' \subset i_* U$, which proves that $\infty \in i_* U$.

\bigskip

Conversely, assume that $\infty \in i_*(U)$, in particular, $i_*(U) \cup X = Y$, hence by locale compactness of $X$: 

\[ Y = \bigcup_{V \ll X} V \cup i_*(U) \]

as $Y$ is compact, there exists $V \ll X$ such that $V \cup i_*(U) = Y$ in particular 

\[X =i^{*}(Y)=i^{*}(V \cup i_*(U) = V \cup i^{*}i_*(U) \leqslant V \cup U \]
which proves $\omega(U)$.

\bigskip

The end of the proof is then a general fact about Artin gluing: consider the natural map $ p : X \coprod \{ \infty \} \rightarrow Y$. This a surjection because $X$ is the open complement of $\{ \infty \}$, hence $\Ocal(Y)$ can be identified with the set of open sublocale $A$ of $ X \coprod \{ \infty \} $ such that $p^{*}p_* A = A$.

An open sublocale of $X \coprod \{ \infty \}$ is exactly a pair $(U \in \Ocal(X), m \in \Ocal(\{ \infty \}) = \Omega)$ and from the first part of the proof one can deduce that $p^{*}p_*(U,m)= (U, m \cap \omega(U))$, which show that $\Ocal(Y)$ is canonically identified with $\Ocal(X^{\infty})$ and there is a unique identification which is compatible to the inclusion of $X$ and $\{\infty \}$.

}

}

\block{In general, a map $f$ from $X^{\infty}$ to any locale $Y$ is the same thing as a map $f_0$ from $X$ to $Y$ and a point $f(\infty) \in Y$ such that for any open sublocale $U \subset Y$ which contains $f(\infty)$ one has $\omega(f^{*}(U))$. Indeed this follows directly from the decomposition of $f^{*} : \Ocal(Y) \rightarrow \Ocal(X^{\infty})$ in the expression of $\Ocal(X^{\infty})$ as a subset of $\Ocal(X) \times \Omega$.

This allows to define:

\Def{Let $\Ccal_{0}(X)$ be the set of functions $f$ from $X$ to the locale $\C$ such that for any positive $\epsilon \in \Q$ one has $\omega(f^{*}(B_{\epsilon} 0))$ where $B_{\epsilon} 0$ denote the ball of radius $\epsilon$ and of center $0$ or, equivalently, the set of functions from $X^{\infty}$ to $\C$ which send $\infty$ to $0$. 
}

}

\block{\Prop{$X^{\infty}$ is completely regular if and only if $X$ is.}

\Dem{If $X^{\infty}$ is completely regular, then any of its sublocales, in particular $X$, is completely regular.

Conversely assume that $X$ is completely regular. Let $A=(U,p) \in \Ocal(X^{\infty})$.
Consider first a $V \ll U$ then there exists $W$ such that $V \triangleleft_{CR} W \ll U \subset X$ and any function on $X$ which is zero outside of $W$ can be extended by $f(\infty)=0$ and hence one has $(V,\bot) \triangleleft_{CR} A$.

Assume now $p$, then one also has $\omega(U)$ hence there exists a $W$ such that $W \ll X$ and $U \cup W =X$. Consider any $W'$ such that $W \ll W' \ll X$, and (as $X$ is completely regular) $f$ a function from $X$ to the locale $[0,1]$ such that $f$ is zero on $W$ and $1$ outside of $W'$. As $W' \ll X$, this function extend to a function from $X^{\infty}$ to $[0,1]$ which satisfies $f^{*}(]0,1]) \subset (U,\top)=(U,p)$ and $f(\infty)=1$. Let $F_{\infty}$ the set of such functions, one can write that:

\[ A= \left( \bigcup_{V \ll U} (V, \bot) \right) \cup  \left( \bigcup_{p,\atop  f \in F_{\infty}} f^{*}(]1/2,1]) \right) \]

which concludes the proof, as, assuming $p$, one has $f^{*}(]1/2,1]) \triangleleft_{CR} A$ for any $f \in F_{\infty}$.

}

}

\block{\label{properfunctoriality}Finally we need to understand the functoriality of the relation between $X$ and $X^{\infty}$:

\Prop{Let $X$ and $Y$ be two regular locally compact locales and $X^{\infty}$ and $Y^{\infty}$ their one point compactifications. Let $f :X \rightarrow Y$ the following conditions are equivalent:

\begin{enumerate}

\item For any $U \in \Ocal(Y)$ such that $U \ll Y$ one has $f^{*}(U) \ll X$

\item There is an extension $f^{\infty}:X^{\infty} \rightarrow Y^{\infty}$ of $f$ such that $(f^{\infty})^{*}(Y)=X$.

\item $f$ is proper (see \cite[C.3.2.5]{sketches}).
\end{enumerate}

Moreover in this case the extension $f^{\infty}$ is unique and this induces a bijection between proper maps from $X$ to $Y$ and maps from $X^{\infty}$ to $Y^{\infty}$ such that $f^{*}(Y)=X$.
}

\Dem{
\begin{itemize}
\item[$1. \Rightarrow 2.$] We define $f^{\infty}$ by:

\[(f^{\infty})^{*}(U,p) = (f^{*}(U),p). \]

Assuming $1.$ if $\omega(U)$ for $U \in \Ocal(Y)$ then there exists a $W$ such that $W \ll Y$ and $W \cup U = Y$, and hence $f^{*}(W) \ll X$ and $f^{*}(U) \cup f^{*}(W)=X$ hence $\omega(f^{*}(U))$. This proves that if $(U,p) \in \Ocal(Y^{\infty})$ , i.e. if $p \Rightarrow \omega(U)$ then one also has $p \Rightarrow \omega(f^{*}(U))$ hence $(f^{*}(U),p) \in \Omega(X^{\infty})$. Moreover as intersections and unions in $\Ocal(X^{\infty})$ are computed componentwise, $f^{\infty}$ is indeed a morphism of locales. As $X$ and $Y$ correspond to the elements $(X, \bot)$ and $(Y,\bot)$ of $\Ocal(X^{\infty})$ and $\Ocal(Y^{\infty})$ one also has $(f^{\infty})^{*}(Y)=X$.

\item[$2. \Rightarrow 3.$] In the situation of $2.$, the map $f$ from $X$ to $Y$ is a pullback of the map $f^{\infty}$ along the open inclusion of $Y$ into $Y^{\infty}$. But any map between two compact regular locales is proper (see \cite[C.3.2.10 (i) and (ii)]{sketches}) and a pullback of a proper map is again a proper map (see \cite[C.3.2.6]{sketches}).

\item[$3. \Rightarrow 1.$] If $f$ is a proper then $f_{*}$ commutes to directed joins. So if we assume that $U \ll Y$, if $X$ is covered by some net $V_i$ then
\[ Y \leqslant f_* \left( \bigcup_i V_i \right) = \bigcup_i f_*(V_i) \]

hence there exists a $j$ such that $U \leqslant f_*(V_j)$ and hence $f^{*}(U) \leqslant V_j$ for some $j$ which concludes the proof of the equivalence.

\end{itemize}

The uniqueness of the extension is immediate because $f^{\infty}$ is defined both on $X$ and on its closed complement, and hence this indeed induce a bijection as stated in the proposition. 
}}

\block{\label{partialfunctoriality}One can be slightly more general:

\Def{Let $X$ and $Y$ be two locally compact regular locales. A partial proper map from $X$ to $Y$ is the data of an open sublocale $\dom(f) \subset X$ and a proper map (denoted $f$) from $\dom(f)$ to $Y$.}

Partial proper maps can be composed (by restricting the domain of definition as much as neccessary) and as a pullback of a proper map is proper the composite of two proper partial maps is again a proper partial map, hence one has a category of proper partial map.

\Prop{The category of pointed compact (completely) regular locales is equivalent to the category of locally compact (completely) regular locales and proper partial maps between them.}

\Dem{The functor are the same as those of proposition \ref{properfunctoriality}, they just apply to a larger category: to a map $f:X^{\infty} \rightarrow Y^{\infty}$ of pointed compact regular locale one associate the partial map $f':X \rightarrow Y$ whose domain is $f^{*}(Y)$, and $f'$ is proper because it is the pullback of $f$ along the inclusion of $Y$ into $Y^{\infty}$ (see the proof of $2. \Rightarrow 3.$ in \ref{properfunctoriality}). Conversely, if $f$ is a partial proper map from $X$ to $Y$ then it extend into a map from $U^{\infty}$ to $Y^{\infty}$ by \ref{properfunctoriality}, composing it to the map $r_U$ of the next lemma yields the desired map from $X^{\infty}$ to $Y^{\infty}$ and these two constructions are clearly inverse of each other.

}

}

\block{\label{mapru}\Lem{Let $X$ be a locally compact regular locale, and $U \subset X$ an open sublocale of $X$ then their exists a (unique) map $r_U : X^{\infty} \rightarrow U^{\infty}$ such that $r_U$ is the identity on $U$ and $(r_U)^{*}(U)=U$.}

\Dem{$r_U$ is defined on $U \subset X^{\infty}$ and on its closed complement (as the constant equal to $\infty$). Hence one has a map  $f_U$ from $U \coprod U^{c}$ to $U^{\infty}$. It is a general fact that the canonical map $U \coprod U^{c} \rightarrow X^{\infty}$ is a surjection of locale (hence corresponds to an injection of frame). So all we have to do to proves that $f_U$ factors into a map $r_U$ on $X^{\infty}$ is to check that for any open sublocale $(V,p) \in \Ocal(X^{\infty})$, the open sublocal $(f_U)^{*}(V,p)$ of $U \coprod U^{c}$ comes from an open sublocale of $U^{\infty}$.

By definition, $(f_U)^{*}(V,p)$ is $V$ on the $U$ part and $s^{*}(p)$ on the $U^{c}$ part (where $s$ is the canonical map $U^{c} \rightarrow {*}$). If we assume $p$ then there exists a $W \ll U$ such that $W \cup V = U$, and as $W \ll U$, there exists a $D \subset X^{\infty}$ such that $W \cap D = \emptyset$ and $U \cup D = X^{\infty}$. We define:

\[ V' = V \cup \coprod_{p \atop D} D \] 

where the coproduct is on the set of $D$ such that $p$ holds and $D$ satisfy the properties just describe. If $f$ denotes the map $U \coprod U^{c} \rightarrow X^{\infty}$ then $p^{*}(V') = (V' \cap U, V' \cap U^{c})$. On one hand $V' \cap U = V$ because as $D \cap W = \emptyset$ and $W \cup V =U$ one has $D \cap U \subset V$, and on the other hand, $U^{c} \cap V' = U^{c} \cap \coprod_{p,D} D $, and as $D  \cup U = X^{\infty}$ , one has $U^{c} \subset D $ hence, $U^{c} \cap V' = s^{*}(p)$, and this concludes the proof.
}

}

\section{Unitarization of $C^{*}$-algebras} 

\label{secUnitarization}

\block{We follow the same definition of $C^{*}$-algebras as for example in \cite{banaschewski2000spectral}. In particular the norm of an element is only assumed to be a upper semi-continuous real number (the definition using rational ball of \cite{banaschewski2000spectral} is equivalent to a norm function with value into the non-negative upper semi-continuous real numbers). Of course contrary to \cite{banaschewski2000spectral}, we do not assume the algebras to be unital.
}

\block{If $C$ is a a $C^{*}$-algebra we define $C^{+}$ as the as the set of couples $(c,z)$ with $c \in C$ and $z \in \C$, We endow $C^{+}$ with the componentwise addition and the multiplication $(c,z)(c',z')=(cc'+cz'+zc',zz')$. It is a unital algebra, with unit $(0,1)$. One also endows $C^{+}$ with the anti-linear involution $(c,z)^{*}=(c^{*},\overline{z})$.

\bigskip

If $(c,z) \in C^{+}$ we define: 

\[ \Vert (c ,z) \Vert = \max( |z|, \sup_{c' \in C_{\leqslant 1} } \Vert c'c+c'z  \Vert )\]

Where $C_{\leqslant 1} $ denotes the set of elements of $C$ of norm $\leqslant 1$, and the supremum is to be considered in the set of upper semicontinuous real numbers. 
}

\blockn{Remark: Classically, it is usual to define the norm on $C^{+}$ to be simply $\sup_{c'} \Vert c'c+c'z  \Vert$. This works perfectly when $C$ is indeed non-unital, but when $C$ is unital this gives a norm $0$ for the element $(1,-1)$ and hence (after taking the quotient by the ideal of norm zero elements) with this definition $C^{+}$ will be isomorphic to $C$ when $C$ is unital. This is not what we want because if $X$ is a compact locale, then its one point compactification $X^{\infty}$ is not $X$ itself but $X \coprod \{ \infty \}$. Classically this difference is harmless but in intuitionist logic the question of being compact/unital might be non decidable and hence it is important to have a uniform treatment on both side.
}

\block{\label{unitarizationLemma1}Before proving that $C^{+}$ is indeed a $C^{*}$-algebra we need a few lemmas which are immediate in classical mathematics but require to be slightly more careful in intuitionist mathematics.

\Lem{Let $x$ be a nonnegative upper semicontinuous real number and $q$ be a nonnegative rational number then if $ x^{2} \leqslant q x $ one has $x \leqslant q $.}

\Dem{ Let $e$ be a rational number such that $x \leqslant q+e$.

One has:

\[ x^{2} \leqslant q^{2} +q e \leqslant (q+\frac{e}{2})^{2} \]

hence $ x \leqslant q+(\frac{e}{2})$

by induction one obtains that for all $k \geqslant 0$ 
\[x \leqslant q+\frac{e}{2^{k}} \]

and hence that $x \leqslant q$ which concludes the proof.

}

}

\block{\label{unitarizationLemma2}\Lem{Let $C$ be a $C^{*}$-algebra and $c\in C$ then:

\[ \Vert c \Vert = \sup_{b \in C_{\leqslant 1}} \Vert bc \Vert = \sup_{b \in C_{\leqslant 1}} \Vert cb \Vert \]
}

\Dem{We start by the first equality. It is immediate that $\Vert bc \Vert \leqslant \Vert c \Vert$ for any $b$ of norm $\leqslant 1$ hence

\[ \sup_{b \in C_{\leqslant 1}} \Vert bc \Vert \leqslant \Vert c \Vert\] 

We only need to prove the reverse inequality. let $q$ be a rational number such that $ \sup_{b \in C_{\leqslant 1}} \Vert bc \Vert < q$, i.e. there exists a $q' <q$ such that for all $b$, $\Vert  bc \Vert <q'$.

Let $\alpha$ be a rational number such that $\Vert c \Vert < \alpha$. One has $\Vert c^{*}/ \alpha \Vert < 1$ hence:

\[\frac{1}{\alpha} \Vert c^{*}c \Vert <q' \]
 
\[ \frac{1}{q'} \Vert c^{*}c \Vert <\alpha, \]

and as this holds for any $\alpha$ such that $\Vert c \Vert < \alpha$ this proves that:

\[ \frac{1}{q'} \Vert c^{*}c \Vert \leqslant \Vert c \Vert \]

Using $\Vert c^{*}c \Vert = \Vert c \Vert^{2}$ one obtains $\Vert c \Vert ^{2} \leqslant q' \Vert c \Vert$ and hence by lemma \ref{unitarizationLemma1} this proves that $\Vert c \Vert \leqslant q' <q$ which concludes the proof of the first equality.

The second equality follows either by exactly the same proof, or by applying the first equality using $\Vert c ^{*} \Vert = \Vert c \Vert$ and $(cb)^{*}=b^{*} c^{*}$.

}

}

\block{\label{unitarizationLemma3}\Lem{For any $x = (c,z) \in C^{+}$ one has:

\[ \sup_{b \in C_{\leqslant 1}} \Vert cb+bz \Vert = \sup_{b \in C_{\leqslant 1}} \Vert bc+bz \Vert, \]

and, $\Vert x \Vert = \Vert x^{*} \Vert$. 

}

\Dem{By lemma \ref{unitarizationLemma2} one has:

\[ \sup_{b \in C_{\leqslant 1}} \Vert cb+bz \Vert  = \sup_{b \in C_{\leqslant 1}} \sup_{b' \in C_{\leqslant 1}} \Vert b'cb+b'bz \Vert \]

But the two supremums can be exchanged and as:

\[ \sup_{b \in C_{\leqslant 1}} \Vert bc+bz \Vert = \sup_{b \in C_{\leqslant 1}} \sup_{b' \in C_{\leqslant 1}} \Vert bcb'+zbb' \Vert \]

This proves the first equality.

The fact that $\Vert x \Vert = \Vert x^{*} \Vert$ follows immediately:

\[ \Vert x^{*} \Vert = \max ( |z|, \sup_{b \in C_{\leqslant 1}} \Vert bc^{*}+b\overline{z} \Vert ) \]

and:

\[ \sup_{b \in C_{\leqslant 1}} \Vert bc^{*}+b\overline{z} \Vert =  \sup_{b \in C_{\leqslant 1}} \Vert c b^{*}+z b^{*} \Vert =  \sup_{b \in C_{\leqslant 1}} \Vert c b+z b \Vert \]

which concludes the proof.

 }

}

\block{\label{unitarizationIsaCstarAlg}\Prop{$C^{+}$ is a $C^{*}$-algebra.}

\Dem{The fact that $\Vert . \Vert$ is a norm of algebra (i.e. such that $\Vert x y \Vert \leqslant \Vert x \Vert \Vert y \Vert $) is easy and exactly as in the classical case. Thanks to the term $ |z|$ in the definition it is an actual norm and not a semi-norm.

$C^{+}$ is complete for this norm because it is complete for the norm $|z|+\Vert c \Vert$ as a product of two Banach spaces, and these two norms are equivalent, indeed, one has in one direction:
\[ \Vert (c,z) \Vert \leqslant | z| + \Vert c \Vert \]
and in the other, one clearly have:
\[ |z | \leqslant \Vert (c,z) \Vert \]
\[ \Vert c \Vert -|z| \leqslant \Vert (c,z) \Vert \]
hence:
 \[ |z | + \Vert c \Vert \leqslant 3 \Vert (c,z)\Vert \]

All we have to do to conclude is to prove the $C^{*}$-equality $\Vert x^{*} x \Vert = \Vert x \Vert^{2}$. Let $x = (c,z)$ and element of $C^{+}$ then:

\[ x^{*} x = (c^{*}c+zc^{*}+\overline{z}c,|z|^{2}) \]

Hence, 
\[ \Vert x^{*} x \Vert =\max (|z|^{2}, \sup_{b \in C_{\leqslant 1}} \Vert c^{*}cb+zc^{*}b+\overline{z}cb \Vert ) \]

as, 
\begin{multline*} \Vert c^{*}cb+zc^{*}b+\overline{z}cb \Vert \geqslant \Vert b^{*}c^{*}cb+zb^{*}c^{*}b + \overline{z}b^{*}cb \Vert \\ = \Vert (cb+z b)^{*} (cb+zb) \Vert =\Vert cb+zb \Vert ^{2} \end{multline*} 

One obtains that $\Vert x^{*} x \Vert \geqslant \Vert x \Vert^{2}$. The other inequality follow from lemma \ref{unitarizationLemma3} together with the fact that $\Vert x^{*} x \Vert \leqslant \Vert x^{*} \Vert \Vert x \Vert$.

}

}

\block{\label{chiinfty}One will identify $C$ with the two sidded ideal of $C^{+}$ of elements of the form $(c,0)$. We will denote by $\chi_{\infty}$ the character of $C^{+}$ defined by $\chi_{\infty}(c,z)=z$.}

\block{\label{unitarizationUniversalProp}
\Prop{Let $C$ be a (possibly non unital) $C^{*}$-algebra and $B$ be a unital $C^{*}$-algebra. Any morphism from $C$ to $B$ extend uniquely into a unital morphism of $C^{*}$-algebra from $C^{+}$ to $B$. }

\Dem{The extension is necessary defined by $f(c,z)= f(c)+z$, it is clearly a morphism of $*$-algebra. It is continuous because:

\[\Vert f(c,z) \Vert \leqslant \Vert f(c) \Vert +|z| \leqslant \Vert f \Vert \Vert c \Vert + |z| \]
and as we observed in the proof of proposition \ref{unitarizationIsaCstarAlg} the norm $\Vert c,z \Vert$ is equivalent to the norm $\Vert c \Vert+ |z|$ this proves that this extension is continuous. }

}

\block{\Prop{Let $X$ be a locally compact completely regular locale, then $\Ccal(X^{\infty}) \simeq (\Ccal_{0}(X))^{+}$ }

\Dem{$\Ccal_{0}(X)$ is identified with the set of functions on $X^{\infty}$ which send $\infty \in X^{\infty}$ to $0$, hence this induce a morphism from $(\Ccal_{0}(X))^{+}$ to $C(X^{\infty})$ by proposition \ref{unitarizationUniversalProp}. A function $f \in C(X^{\infty})$ can be written in a unique way $h+c$ with $h \in \Ccal_{0}(X)$ and $c$ a constant: $c$ has to be $f(\infty)$ and $h = f - f(\infty)$ hence the map from $(\Ccal_{0}(X))^{+}$ to $C(X^{\infty})$ is a bijection. One easily check that it is isometric, either by general theorems on $C^{*}$-algebras (which should of course be proved constructively first) or directly: If $f$ is a function on $X^{\infty}$ then $\Vert f \Vert <q$ if and only if both $f(\infty)<q$ (which imply that $f < q$ on some neighbourhood of $\infty$) and for each $U \triangleleft_{CR} X$ the function $f$ is strictly smaller than a $q'<q$ on $U$, which is equivalent to the fact that $\Vert fh \Vert <q'<q$ for every $h \in \Ccal_0(X)$ of norm $\leqslant 1$. And these two conditions are equivalent to the fact that $\Vert (f-f(\infty),f(\infty) ) \Vert <q$.
}
}

\block{We conclude this section by discussing the compatibility of unitarization to pullback along geometric morphisms.

If $\Ccal$ is a $C^{*}$-algebra in a topos $\Ecal$ and $f:\Tcal \rightarrow \Ecal$ is a geometric morphism, then $f^{*}(\Ccal)$ is in general not a $C^{*}$-algebra: it still satisfies all the ``algebraic" axioms. But it might not be complete and separated (in the sense that $\Vert x \Vert =0 \Rightarrow x=0$). Fortunately, the separated completion of $f^{*}\Ccal$ is complete and separated and hence is a $C^{*}$-algebra which we denote\footnote{This also corresponds to the pullback of the localic completion of $\Ccal$, hence this is essentially compatible with the notation for pullback of locales} by $f^{\sharp}(\Ccal)$.
\Prop{For any $C^{*}$-algebra $\Ccal$ one has a natural isomorphism:

\[ f^{\sharp}( \Ccal^{+}) \simeq f^{\sharp}(\Ccal)^{+} \]

}
\Dem{In $\Tcal$, the canonical morphisms of pre-$C^{*}$-algebras $f^{*} \Ccal \rightarrow f^{\sharp}(\Ccal)$ and $f^{*}\C \rightarrow \C$ extend into a map $f^{*}\C \times f^{*}\Ccal \rightarrow f^{\sharp}(\Ccal)^{+}$. One can check that the semi-norm and the pre-$C^{*}$-algebra structure induced on $\C \times f^{*} \Ccal$ by this map are exactly those of $f^{*}(\Ccal^{+})$, hence $f^{\sharp}(\Ccal^{+})$ is exactly the closure of $f^{*}(\Ccal^{+})$ in $f^{\sharp}(\Ccal)^{+}$. But $f^{*}(\Ccal^{+})$ is clearly dense (because each component is dense) hence this concludes the proof. }

}

\section{The non-unital Gelfand duality}

\block{\Def{If $\Ccal$ is a $C^{*}$-algebra we denote by $\spect \Ccal$ the spectrum of the unital $C^{*}$-algebra $\Ccal^{+}$ and by $\spec \Ccal$ the locally compact completely regular locale obtained by removing the point $\infty$ of $\spect \Ccal$.}

Of course by the uniqueness property in theorem \ref{OnepointCptmain}, $\spect \Ccal$ is the one point compactification of $\spec \Ccal$. Also, if $\Ccal$ is unital then $\Ccal^{+}$ is isomorphic to $\Ccal \times \C$ hence $\spect \Ccal$ is isomorphic to $\spec \Ccal \coprod \{ \infty \}$ and the two definitions of $\spec \Ccal$ (by considering $\Ccal$ as a unital or general $C^{*}$-algebra) agree and there is no possible confusion.
}

\blockn{At this point, the following theorem is immediate:}

\block{\label{mainResult}\Th{The category of commutative $C^{*}$-algebras and arbitrary morphisms between them is anti-equivalent to the category of locally compact completely regular locales and partial proper maps between them. The equivalence is given on object by the constructions $\Ccal_0$ and $\spec$.}

\Dem{The process of unitarization produce an equivalence between the category of commutative $C^{*}$-algebra and arbitrary morphism, and the category of unital $C^{*}$-algebras endowed with a character $\chi_{\infty}$ and unital morphism compatible to the character. Applying the Gelfand duality for unital $C^{*}$-algebra this category is in turn anti-equivalent to the category of pointed compact completely regular locales, which by proposition \ref{partialfunctoriality} is equivalent to the category of locally compact completely regular locale and partial proper map between them. Under these composed equivalences, a commutative $C^{*}$-algebra $\Ccal$ is associated to the spectrum of $\Ccal^{+}$ minus the point at infinity, i.e. exactly $\spec \Ccal$ and a locally compact completely regular locale $X$ is associated to the algebra of functions on $X^{\infty}$ which vanish at $\infty$, which is $\Ccal_0(X)$.
}

In the rest of this section, we will give interpretation of $\spec$ and $\spect$ in term of classifying space of characters (proposition \ref{spec_class_char}), we will show that the open sublocales of $\spec \Ccal$ correspond to the closed ideals of $\Ccal$ (theorem \ref{openeqideal}) and that total proper maps of locales correspond to non-degenerate morphisms of $C^{*}$-algebras (theorem \ref{non-degen=propermap}).

}

\block{We recall that when $\Ccal$ is a unital commutative $C^{*}$-algera, then $\spec \Ccal$ denotes the classyfing space of the theory of characters of $\Ccal$. A precise geometric formulation of this theory can be found in \cite{banaschewski2000spectral} or in \cite{coquand2009constructive}, but this can also be interpreted as the fact that for any locale $Y$ and $p:Y \rightarrow \{*\}$ the canonical map, functions from $Y$ to $\spec \Ccal$ correspond to morphisms from $p^{\sharp}(\Ccal)$ (or equivalently from $p^{*}(\Ccal)$) to $\C$ internally in $\sh(Y)$.
}

\block{\label{spec_class_char}
\Prop{The locale $\spect \Ccal$ classifies ``nonunital characters" of $\Ccal$, i.e. possibly nonunital morphism of $C^{*}$-algebras from $\Ccal$ to $\C$.

The locale $\spec \Ccal$ classifies ``nonzero characters" of $\Ccal$, i.e. characters which satisfy the additional axiom $\exists c \in \Ccal, |\chi(c)|>0$
}
 
\Dem{The important observation is that the process of unitarization of $C^{*}$-algebras commute with pullback along geometric morphisms. Hence points of $\spect \Ccal$ over any locale $\Lcal$ (with $p$ its canonical morphism to the point) are the characters of $p^{\sharp}(\Ccal)^{+}$ which are in bijection with the non unital morphisms from $p^{\sharp}(\Ccal)$ to $\C$ which proves the first part of the result. 

For the second part, let us denote by $D(f)$ the (biggest) open sublocale of $\spect \Ccal$ on which $|f|>0$, where $f$ is an element of $\Ccal$. Using the complete regularity of $\spect \Ccal$ it appears that $\spec \Ccal$ is the union of the $D(f)$ for $f \in \Ccal$. For each $f \in \Ccal$, the open sublocale $D(f)$ classifies characters of $\Ccal$ such that $|\chi(f)|>0$. Hence points of $\spec \Ccal$ are the characters such that $\exists c \in p^{*}(\Ccal) $ with $|\chi(c)|>0$.

The formulation ``$\exists c \in \Ccal$" in the statement of the proposition is unambiguous because as $p^{*}(\Ccal)$ is dense in $p^{\sharp}(\Ccal)$ it is equivalent to says that $\exists c \in p^{*}(\Ccal), |\chi(c)|>0$ and that  $\exists c \in p^{\sharp}(\Ccal), |\chi(c)|>0$ hence this concludes the proof.
}

}

\block{\label{lemmaC0U}\Lem{Let $X$ be a locally compact completely regular locale and $U \subset X$ an open sublocale, then the restriction to $U$ of functions in $\Ccal_0(X)$ which vanish outside of $U$ are exactly the functions in $\Ccal_0(U)$.}

\Dem{Let $f$ be a function in $\Ccal_0(X)$ which vanish outside of $U$, we will show that the restriction of $f$ to $U$ is in $\Ccal_0(U)$.
Let $\epsilon$ be any positive rational number, let $V$ be the open sublocale of $X$ on which $|f|>\epsilon$. As $f \in \Ccal_0(X)$, $V\ll X$ and as $f$ vanish outside of $U$ on has $V \triangleleft_{CR} U$ and in particular $V \triangleleft U$. As mentioned in the preliminaries, this two properties together imply $V \ll U$. This being true for any $\epsilon$ this proves that $f \in \Ccal_0(U)$.

Conversely, assume that $f\in \Ccal_0(U)$, then because of the map $r_U: X^{\infty} \rightarrow U^{\infty}$ constructed in \ref{mapru} one can define a map $f\circ r_U$ on $X ^{\infty}$ which vanish at infinity and outside of $U$ and which coincide with $f$ on $U$.
}
}

\block{\label{lemmaextcharac}\Lem{Let $I \subset \Ccal$ be an ideal of a commutative $C^{*}$-algebra. and let $\chi$ be a non-zero character of $I$ (in the sense that $\exists i \in I, |\chi(i)|>0 $). Then $\chi$ admit a unique extension as a character of $\Ccal$.}

\Dem{The proof is exactly as in the classical case\footnote{Except that we need to be more careful on the ``non-zero" hypothesis which is unessential in the classical case.}: any extension of $\chi$ to $\Ccal$ has to satisfy $\chi(c)=\chi(c i)/\chi(i)$ for any $i \in I$ such that $|\chi(i)|>0$ hence the extension is unique. Conversely, if one has $f,g$ two elements of $I$ such that $|\chi(f)|>0$ and $|\chi(g)|>0$, and $c$ any element of $\Ccal$ then, as:

\[ \chi(gfc) = \chi(g) \chi(fc) = \chi(f) \chi(gc) \] 
one has:
\[ \frac{\chi(fc)}{\chi(f)} = \frac{\chi(gc)}{\chi(g)}. \]

This proves that we can define $\chi(c)=\chi(fc)/\chi(f)$ for any $f \in I$ such that $\chi(f)$ is invertible and as $\chi(c)\chi(c')=\chi(fc)\chi(fc')/\chi(f)^{2}= \chi(f^{2} cc')/\chi(f^{2}) = \chi(cc')$, the extension of $\chi$ is a character of $\Ccal$.
}
}

\block{\label{openeqideal}\Th{The constructions $\Ccal_0$ and $\spec$induce for any commutative $C^{*}$-algebra $\Ccal$ an order preserving bijection between the open sublocales of $\spec \Ccal$ and the closed ideals of $\Ccal$.
\smallskip

Moreover, if $f:\Ccal \rightarrow \Ccal'$ is a morphism between commutative $C^{*}$-algebra the pull-back of open sublocales along the corresponding (partial) continuous map corresponds under this bijection to the map which send an ideal $I$ of $\Ccal$ to the closure of the ideal spammed by $f(I)$. }

\Dem{Because $\Ccal$ is an ideal of $\Ccal^{+}$ it suffices to prove these results for unital algebras and a unital morphism.

If $U \subset \spec \Ccal$ is an open sublocale, then by lemma \ref{lemmaC0U} one can identify $\Ccal_0(U)$ with an ideal of $\Ccal$ whose spectrum is $U$. Conversely, if $I \subset \Ccal$ is an ideal then an application of lemma \ref{lemmaextcharac} internally in $\spec I$ give rise to a map from $\spec I$ to $\spec \Ccal$ and a map from an arbitrary locale $\Lcal$ to $\spec \Ccal$ factor into $\spec I$ if and only the corresponding character of $\Ccal$ in the logic of $\Lcal$ satisfies $\exists i\in I, |\chi(i)|>0$ (this is again an application of \ref{lemmaextcharac} internally to $\Lcal$). Hence $\spec I$ is identified precisely with the open sublocales of $\Ccal$ defined by $\bigcup_{i \in I} D(i)$. And as $\Ccal_0(\spec I) = I$ this proves that this two constructions are inverse of each other, and they clearly preserve the order.

For the second part, if $ f:\Ccal \rightarrow \Ccal'$ is a unital morphism of $C^{*}$-algebra, $g$ the corresponding continuous map $\spec \Ccal' \rightarrow \spec \Ccal$, and if $I \subset \Ccal$ and $I' \subset \Ccal'$ are two closed ideals, then $f(I) \subset I'$ if for any function $h$ on $\spec \Ccal$ which vanish outside of $\spec I$ its composite with $g$ vanish outside of $\spec I'$. This will be the case if and only if $g^{*}(\spec I) \subset \spec I'$.  Hence the ideal corresponding to $g^{*}(\spec I)$ is indeed the smallest closed ideal containing $f(I)$.
}
}

\block{\label{non-degen=propermap}The following theorem is in fact just a corollary of theorem \ref{openeqideal}. We recall that a morphism of $C^{*}$-algebra is said to be non-degenerate it its image spam a dense ideal. A morphism between unital $C^{*}$-algebra is non-degenerate if and only if it is unital.

\Th{The equivalence of category of theorem \ref{mainResult} restrict to a (contravariant) equivalence of category between the category of commutative $C^{*}$-algebras and non degenerate morphisms and the category of locally compact completely regular locales and proper maps between them.}

\Dem{A morphism $f:\Ccal \rightarrow \Ccal'$ will corresponds to a total map on the spectrum if and only if the corresponding partial proper map $g$ satisfy $g^{*}(\spec \Ccal)= \spec \Ccal'$, i.e., applying the previous theorem, if and only if $f$ is non-degenerate.}

}

\section{Local positivity and continuity of the norm}

\block{We recall that a locale $X$ is said to be positive if when $X = \bigcup_{i \in I} U_i$ then $\exists i \in I$. This is a ``positive" way of saying that $X$ is non-zero. We will say that a locale $X$ is locally positive if every open sublocale of $X$ can be written as a union of positive open sublocales. Assuming the law of excluded middle, a locale is positive if and only if it is non-zero and any locale is locally positive, but in an intuitionist framework locale positivity is an extremly important properties: A locale $X$ is locally positive if and only if the map from $X$ to the terminal locale is an open map (see \cite[C3.1.17]{sketches}) for this reason locally positive locale are aslo called open locale (but this cause a confusion with open sublocales) or sometimes overt.}

\block{\label{continuity=openess}\Th{Let $\Ccal$ be a commutative $C^{*}$-algebra, then the following conditions are equivalent:

\begin{itemize}
\item For any $c \in \Ccal$, the norm of $c$ is a continuous real number.
\item There is a dense family of elements of $\Ccal$ whose norms are continuous real numbers.
\item $\spec \Ccal$ is locally positive. (i.e. is open or overt).
\end{itemize}
}

It appears that this result was already known for unital algebras and due to T.Coquand in \cite[section 5]{coquand2005stone}. We do needed the result for non-unital algebras in \cite{henry2015toward}, but one could also deduce the non-unital case from the unital one using the unitarization process developed in section \ref{secUnitarization}. This being said, we were not aware of Coquand's paper at the time the first version of this paper has been written, and as the following proof is more complete than the original one we decided to leave it here.

\Dem{The first two conditions are clearly equivalent because a semi-continuous real number which can approximated arbitrarily closed by continuous real numbers is also continuous.

Assume first that the first two conditions hold.

We recall that if $f \in \Ccal$ then $D(f)$ denotes the largest open sublocale of $\spec \Ccal$ on which $|f|>0$. Let also $p$ denotes the canonical map from $\spec \Ccal$ to the terminal locale.

We will first show that:

\[ D(f) \subset p^{*}(``\Vert f \Vert>0"). \] 
Indeed, in the logic of $\spec \Ccal$, $D(f)$ is the proposition $\chi(|f|)>0$, which imply that $\exists \epsilon \chi(|f|)>\epsilon >0$. But, as $\Vert f \Vert $ is continuous, one has (still internally in $\spec \Ccal$) $\Vert f \Vert <\epsilon$ or $\Vert f \Vert >0$. $\Vert f \Vert <\epsilon$ is in contradiction with $\chi(|f|)>\epsilon$, hence $\Vert f \Vert >0$.

The $(D(f))_{f \in \Ccal}$ form a basis of the topology of $\spec \Ccal$, hence as:
\[D(f) = D(f) \cap p^{*}(``\Vert f \Vert>0") = \bigcup_{\Vert f \Vert >0} D(f) \]
the $D(f)$ for $\Vert f \Vert >0$ also form a basis of the topology of $\spec \Ccal$. We will now prove that the $D(f)$ for $\Vert f \Vert >0$ are positive and this will conclude the proof of this implication.

If $\Vert f \Vert >0$ there exists a rational $\epsilon>0$ such that $\Vert f \Vert >\epsilon$. Let $U_i$ be a familly of open sublocales of $D(f)$ such that:

\[ D(f) = \bigcup_{i \in I} U_i \]

Let $W$ be the open sublocale on which $|f|$ is greater than $\epsilon/2$. One has $W \triangleleft_{CR} D(f)$ by definition and $W \ll \spec \Ccal$ because $f \in \Ccal_0(\Ccal)$, hence, as mentioned in the preliminaries, one has $W \ll D(f)$, hence there exists a finite subset $J \subset I$ such that:

\[ W \subset \bigcup_{j  \in J} U_j \]

As $J$ is finite, it is either empty or inhabited, but if $J$ is empty then $W$ is empty hence $f$ is smaller than $\epsilon/2$ everywhere on $\spec \Ccal$ and hence $\Vert f \Vert < \epsilon$ which yields a contradiction. This shows that $J$ is inhabited and hence that $I$ is inhabited, which concludes the proof of the first implication.

\bigskip

We now assume that $\spec \Ccal$ is a locally positive locale. For any $h \in \Ccal$, we denote $(|h|>q)$ the biggest open sublocale of $\spec \Ccal$ on which $|h|>q$ holds (where $h$ is seen as a function on $\spec \Ccal$). We fix an element $h \in \Ccal$, and we will prove that $\Vert h \Vert$ is a continuous real number. We define:

\[ L = \{q \in \Q | q<0 \text{ or } (|h|>q) \text{ is positive } \} \]

one has:

\begin{itemize}

\item If $q\in L$ and $q' < q$ then $q' \in L$
\item $L$ is inhabited (it contains all the negative rational numbers).
\item if $q \in L$ then there exists $q' \in L$ such that $q<q'$, indeed, if $q<0$ it is clear, and if $(|h|>q)$ is positive then it is the union for $q'>q$ of the $(|h|>q')$, as $\spec \Ccal$ is assumed to be locally positive $(|h|>q)$ is also the union of the $(|h|>q')$ which are positive and hence there exists a $q'$ such that $(|h|>q')$ is positive.
\end{itemize}

This shows that $L$ is a lower semi-continuous real number. We will show that $(L, \Vert h \Vert )$ form a continuous real number, which means that $\Vert h \Vert$ is a continuous real number.

\begin{itemize}

\item Let $q$ such that $q \in L$ and $\Vert h \Vert <q$, this means that $|h|$ is both smaller than $q$ everywhere and bigger than $q$ on some positive sublocale which is impossible. Hence $L \cap \Vert h \Vert = \emptyset$.

\item Let $q<q'$ be two rational numbers. Internally in $\spec \Ccal$ one has $|h|<q'$ or $q<|h|$. Hence $\spec \Ccal$ is the union of the open sublocales $(|h|<q')$ and $(q<|h|)$, moreover as $h\in \Ccal_0(\spec \Ccal)$ one has $(q<|h|) \ll \spec \Ccal$.
By locale positivity of $\spec \Ccal$, the open sublocale $(q<|h|)$ can be written as a union positive open sublocales $(u_i)$ for $i \in I$. In particular:

\[\spec \Ccal = (|h|<q') \cup \bigcup_{i \in I} u_i \]

hence, as $(q<|h|) \ll \spec \Ccal$, there exists a finite subset $J \subset I$ such that:
\[ (q<|h|) \subset (|h|<q') \cup \bigcup_{j \in J} u_j \]

As $J$ is finite, it is either empty or inhabited. If $J$ is empty, then $(q<|h|) \subset (|h|<q')$ hence $(|h|<q') = ([h|<q') \cup (q<|h|) = \spec \Ccal$ hence $\Vert h \Vert < q'$. On the other hand, if $J$ is inhabited then $(q < |h|)$ contains a positive open sublocale, hence it is positive and hence $q \in L$. This proves that either $\Vert h \Vert<q'$ or $q \in L$.

\end{itemize}

This two conditions together show that $(L,\Vert h \Vert)$ form a continuous real number and this concludes the proof.

}

}

\section{Extension of the results to localic $C^{*}$-algebras}

\label{secLocalic}

\block{In \cite{henry2014localic} we have defined a notion of ``localic $C^{*}$-algebras" and proved (as previously conjectured by C.J.Mulvey and B.Banachewski in \cite{banaschewski2006globalisation}) that the (constructive) Gelfand duality can be extended into a duality between compact regular locales and localic commutative unital $C^{*}$-algebras. The goal of this last section is to explain how the methods developed in \cite{henry2014localic} allow to extend the results of the present paper to the localic framework (we have in mind theorems \ref{mainResult}, \ref{openeqideal} and \ref{continuity=openess} and proposition \ref{spec_class_char}). In particular, this section is not meant to be read independently of \cite{henry2014localic}.}

\block{Let us start with the construction of the spectrum of a localic $C^{*}$-algebra.

\Prop{If $\Ccal$ is a (possibly non-unital) commutative $C^{*}$-locale then there exist locales $\spect \Ccal$ and $\spec \Ccal$ such that $\spect(\Ccal)$ classifies the morphisms $\chi : \Ccal \rightarrow \C$ of $C^{*}$-locales, and $\spec \Ccal$ classifies those which satisfy additionally ``$\chi^{-1}(\C-\{0\})$ is positive". Moreover $\spec \Ccal$ is a locally compact regular locale and  $\spect(\Ccal)$ is its one point compactification.
Finally, these constructions are compatible with pullback along geometric morphisms. }

\Dem{In section $3.5$ of \cite{henry2014localic} we proved that there is a classifying space for metric maps from $\Ccal$ to $\C$ denoted $[ \Ccal,\C]_1$, one can then construct $\spect \Ccal$ as a sublocale of $[\Ccal,\C]_1$ using the same kind of co-equalizer as in $4.2.3$ of \cite{henry2014localic}. Moreover ``$\chi^{-1}(\C-\{0\})$ is positive " is an open subspace of $[\Ccal,\C]_1$ (it is even one of the basic open subspace) hence $\spec \Ccal$ will be an open subspace of $\spect \Ccal$.

The compatibility with pullback along geometric morphisms follows immediately from this definition as classifying space.

For the rest of the proposition we can use descent theory: from proposition $2.3.17$ of \cite{henry2014localic} there exists a locale $\Lcal$, with $p:\Lcal \rightarrow \{*\}$ the canonical map, such that $p^{\sharp}(\Ccal)$ is weakly spatial, hence is the localic completion of an ordinary $C^{*}$-algebra. In particular, as character of a $C^{*}$-algebra and of its localic completion are the same, $p^{\sharp}(\spec \Ccal)$ and $p^{\sharp}(\spect \Ccal)$ are the spectrums of an ordinary $C^{*}$-algebra, hence they are respectively locally compact completely regular and compact completely regular and the second is isomorphic to the one point compactification of the first. (Complete) regularity alone is not a property that descend well along open surjections, but it is proved in \cite[Lemma C3.2.10]{sketches} that for compact locale regularity is equivalent to begin Hausdorff and $\cite[C5.1.7]{sketches}$ prove that as $p^{\sharp}(\spect \Ccal)$ is compact and separated, $\spect \Ccal$ is also compact and separated, hence compact regular. As $\spec \Ccal$ is an open subspace of $\spect \Ccal$ it is locally compact and regular. Finally, as one point compactification is also compatible with pullback along geometric morphism the isomorphism between $p^{\sharp}(\spect \Ccal)$ and the one point compatification of $p^{\sharp}(\spec \Ccal)$ descend into an isomorphism between $\spect \Ccal$ and the one point compactification of $\spec \Ccal$ (which of course is compatible with the natural inclusion of $\spec \Ccal$).

}

}

\block{\Th{There is an anti-equivalence of categories between the category of commutative (possibly non-unital) $C^{*}$-locales and the category of locally compact regular locales and partial proper maps between them. }

\Dem{The proof given in \cite[4.2.5]{henry2014localic} that the ordinary Gelfand duality extend to the localic Gelfand duality applies to the non-unital case almost without any change: If $X$ is a locally compact regular locale one can define the $C^{*}$-locale $\Ccal_0(X)$ as the kernel of the evaluation at infinity on the $C^{*}$-locale $\Ccal(X^{\infty})$. $\Ccal_0(X)$ is a $C^{*}$-locale: the only non trivial point to check is that it is locally positive but it follows from the fact that the map $\Ccal(X^{\infty}) \rightarrow \Ccal_0(X)$ which send $f$ to $f-f(\infty)$ is a surjection.  There is a canonical map from $X$ to $\spec \Ccal_{0}(X)$. By \cite[2.3.17 and 2.6]{henry2014localic} there exists a positive locally positive locale $\Lcal$ such that $p^{\sharp}(X^{\infty})$ is completely regular in the internal logic of $\Lcal$ (and hence also $p^{\sharp}(X)$), in particular the canonical map $p^{\sharp}(X) \rightarrow p^{\sharp}(\spec \Ccal_{0}(X)) \simeq \spec \Ccal_{0}(p^{\sharp}(X))$ is an isomorphism because of the ordinary non-unital gelfand duality applied to the $C^{*}$-algebra of points of $\Ccal_0(X)$, and because open surjections are effective descent morphisms (and $p$ is an open surjection) this imply that $X \simeq \spec \Ccal_{0}(X)$.

The exact same argument also show that for any commutative $\Ccal^{*}$-locale $\Ccal$ the canonical map $\Ccal \rightarrow \Ccal_0(\spec \Ccal)$ is an isomorphism, and the correspondence of morphisms is also obtains exactly in the same way (because partial maps descend well : first apply descent to their domain and then to the map itself).
}  }

\block{\Th{Let $\Ccal$ be a commutative $C^{*}$-locale, then there is an order preserving bijection between open sublocales of $\spec \Ccal$ and locally positive fiberwise closed ideals of $\Ccal$.} 

One might be surprised to obtain ``fiberwise" closed ideals, and not just closed ideals in this duality. But there is no reason to be, one should just notice that the usual notion of closeness ans density that we use in the constructive theory of Banach space does not correspond to closeness and density but to fiberwise closeness and fiberwise density.

Indeed, a point $x$ is in the closure of a subset $S$ of a Banach space $B$ if for all $\epsilon >0$ there exists a $s \in S$ such that $\Vert x - s \Vert < \epsilon$, i.e. such that for all neighbourhood $V$ of $x$, there is a point in $V \cap S$, i.e. in localic terms, $V \cap S$ is positive which exactly says that $x$ is in the fiberwise closure of $S$.

\Dem{Let $\Ccal$ be a commutative $C^{*}$-locale, and let $I$ be a locally positive fiberwise closed ideal of $\Ccal$. In particular $I$ is a $C^{*}$-locale and hence has a spectrum $ U = \spec I$. Let $\Lcal$ be a positive locally positive locale such that $p^{\sharp}(I)$ and $p^{\sharp}(\Ccal)$ are weakly spatial, by \ref{openeqideal} $p^{\sharp}(U)$ identify with a open sublocale of $p^{\sharp}(\spec  \Ccal)$. This open injection is canonical hence compatible to the descent data and hence comes from a map from $U$ to $\spec \Ccal$ which also has to be an open injection (because its pullback along $p$ is an open injection). Conversely, if $U$ is an open sublocale of $\spec \Ccal$ then $C_{0}(U)$ is a $C^{*}$-locale and identify by the same descent argument with a locally positive fiberwise closed ideal of $\Ccal$, and these two construction are inverse of each other essentially because of the localic gelfand duality we just proved. }
}

\blockn{And finally:}

\block{\Th{Let $\Ccal$ be a commutative $C^{*}$-locale, then the following conditions are equivalent:

\begin{itemize}
\item $\spec \Ccal$ is locally positive.
\item the norm map from $\Ccal$ to the locale of upper semi-continuous real number factors into the natural map from the locale of continuous real number to the locale of upper-semi-continuous number.
\end{itemize}
}
\Dem{Let $\Lcal$ be a positive locally positive locale and $p$ the canonical map $p:\Lcal \rightarrow \{* \}$ such that $p^{\sharp}(\Ccal)$ is weakly spatial.

\bigskip

Assume the first condition, then in the logic of $\Lcal$, the locale $\spec p^{\sharp}\Ccal \simeq p^{\sharp} \spec \Ccal$ is still locally positive, hence by theorem \ref{continuity=openess} each point of $p^{\sharp}\Ccal$ has a continuous norm. The subalgebra of points is fiberwise dense and endowed with a map to the locale of continuous real number which factors the norm map. This map from the locale of points to the locale of real numbers is clearly a metric map and hence extend by completion to a map from the all of $p^{\sharp}\Ccal$ to the locale $\R$ which also factors the norm (by \cite[3.3 and 3.2.5]{henry2014localic}). By uniqueness of such a factorisation, it is compatible with the descent datas on $p^{\sharp}\Ccal$ and $p^{\sharp}\R$ and hence induces a map from $\Ccal$ to $\R$ which also factor the norm and this concludes the proof of the first implication.

\bigskip

We now assume the second condition. This factorisation of the norm implies that every point of $p^{\sharp}(\Ccal)$ has a continuous norm, hence $\spec p^{\sharp}(\Ccal) \simeq p^{\sharp}(\spec \Ccal)$ is locally positive by \ref{continuity=openess} and hence $\spec \Ccal$ is also locally positive by \cite[C5.1.7]{sketches}.
} }

\bibliography{Biblio}{}
\bibliographystyle{plain}

\end{document}